\documentclass[reqno,11pt]{article}
\usepackage{mathrsfs}
\usepackage{amssymb}
\usepackage{amsthm}
\usepackage{amsmath}
\usepackage{amsfonts}
\usepackage{enumerate}
\usepackage{bm}

\usepackage{graphicx}
\usepackage{setspace,color,wrapfig}
\usepackage[colorlinks,linkcolor=blue,citecolor=blue]{hyperref}

\numberwithin{equation}{section}
\newtheorem{theorem}{Theorem}[section]
\newtheorem{lemma}[theorem]{Lemma}

\theoremstyle{definition}
\newtheorem{definition}[theorem]{Definition}

\theoremstyle{remark}
\newtheorem{remark}[theorem]{Remark}

\begin{document}
\title{Two dimensional subsonic and subsonic-sonic spiral flows outside a porous body }
\author{Shangkun Weng\thanks{School of mathematics and statistics, Wuhan University, Wuhan, Hubei Province, 430072, People's Republic of China. Email: skweng@whu.edu.cn}\and Zihao Zhang\thanks{School of mathematics and statistics, Wuhan University, Wuhan, Hubei Province, 430072, People's Republic of China. Email: zhangzihao@whu.edu.cn}}
\date{}
\maketitle
\newcommand{\de}{{\mathrm{d}}}
\def\div{{\rm div\,}}
\def\curl{{\rm curl\,}}
\newcommand{\ro}{{\rm rot}}
\newcommand{\sr}{{\rm supp}}
\newcommand{\sa}{{\rm sup}}
\newcommand{\va}{{\varphi}}
\newcommand{\me}{\mathcal{M}}
\newcommand{\ml}{\mathcal{V}}
\newcommand{\mi}{\mathcal{N}}
\newcommand{\md}{\mathcal{D}}
\newcommand{\mg}{\mathcal{G}}
\newcommand{\mh}{\mathcal{H}}
\newcommand{\mf}{\mathcal{F}}
\newcommand{\n}{\nabla}
 \newcommand{\q}{{\rm R}}
\newcommand{\p}{{\partial}}

\begin{abstract}
In this paper, we investigate two dimensional subsonic and subsonic-sonic spiral flows outside a porous body. The existence and uniqueness of the subsonic spiral flows are obtained via variational formulation. The optimal decay rate at far field is also derived by the Kelvin's transformation and some elliptic estimates. By extracting spiral subsonic solutions as the approximate sequences, we obtain the spiral subsonic-sonic limit solution. The main ingredients of our analysis  are methods of calculus of variations, the theory of second-order quasilinear equations and the compactness framework.
\end{abstract}

\begin{center}
\begin{minipage}{5.5in}
Mathematics Subject Classifications 2010: Primary 35B40, 35Q31; Secondary 35J25, 76N15.\\
Key words:  subsonic spiral flows, Euler equations, subsonic-sonic limit, a porous body.
\end{minipage}
\end{center}

\section{Introduction and main results}\noindent
\par In this paper, we are concerned with two-dimensional subsonic and subsonic-sonic spiral flows outside a porous body $ \md(\Gamma) $, which are governed by the following two-dimensional Euler system:
\begin{equation}\label{1-1}
\begin{cases}
\p_{x_1}(\rho u_1)+\p_{x_2}(\rho u_2)=0,\\
\p_{x_1}(\rho u_1^2)+\p_{x_2}(\rho u_1u_2)+\p_{x_1}p=0,\\
\p_{x_1}(\rho u_1u_2)+\p_{x_2}(\rho u_2^2)+\p_{x_2}p=0,\\
\end{cases}
\end{equation}
where $ \textbf{u} = (u_1, u_2) $ is the velocity field, $ \rho $ is the density, and $ p $ is pressure. Here we only consider the polytropic gas, therefore $ p=A\rho^{\gamma} $, where  $ A $ is a positive constant and $ \gamma $ is the adiabatic constant with $ \gamma> 1 $.
\par Suppose that the flow is also irrotational, i.e.
\begin{equation}\label{1-2}
\p_{x_2}u_1-\p_{x_1}u_2=0.
\end{equation}
Then it follows from \eqref{1-1} and \eqref{1-2} that the flow satisfies  the Bernoulli's law
\begin{equation}\label{1-3}
\frac{q^2}{2}+h(\rho)=C_0,
\end{equation}
where  $ h(\rho) $ is the enthalpy satisfying $ h^{\prime}(\rho)=\frac{p^{\prime}(\rho)}{\rho} $, $ q=\sqrt{u_1^2+u_2^2} $ is the flow speed and $ C_0 $ is a constant depending on the flow.   We  normalize the flow as in \cite{CF48}, such that $ p(\rho)=\rho^{\gamma}/\gamma $ is the pressure for the polytropic gas, and the Bernoulli's law \eqref{1-3} reduces to
\begin{equation}\label{1-4}
\frac{q^2}{2}+\int_0^{\rho}\frac{p^{\prime}(\rho)}{\rho}\mathrm{d}\rho
=\frac{\gamma+1}{2(\gamma-1)}.
\end{equation}
So \eqref{1-4} yields a representation of the density
 \begin{equation}\label{1-5}
 \rho=g(q^2)=\left(\frac{\gamma+1-(\gamma-1)q^2}{2}\right)^{1/(\gamma-1)}.
 \end{equation}
 The sound speed $ c $  is defined as $ c^2=p^{\prime}(\rho)$. At the sonic point $ q = c $, \eqref{1-5} implies $  q^2=1 $. We define the critical speed $ q_{cr} $ as $  q_{cr}=1 $. Thus the flow is subsonic when $ q<1 $, sonic when $ q=1 $ and supersonic when $ q>1 $.
 \par By the density equation in \eqref{1-1}, one can introduce the stream function $ \psi $ as follows
\begin{equation}\label{1-6}
\p_{x_1}\psi =-\rho u_2,\quad \p_{x_2}\psi=\rho u_1.
\end{equation}
 Obviously, $ |\n\psi|=\rho q $. Then it follows from the normalized  Bernoulli's law \eqref{1-4} that $ \rho q $ is a nonnegative function of $ q $,  which is
increasing for $ q \in (0, 1 ) $ and decreasing for $ q \geq 1 $, and vanishes at $ q=0 $. So $ \rho q $ attains its maximum
at $ q=1 $.  Therefore $ \rho $ is a two-valued function of $ |\n\psi|^2 $. Subsonic flows correspond to the branch where $\rho > 1 $  if $ |\n\psi|^2 \in [0, 1) $. Set
 \begin{equation}\label{1-7}
  \rho=H(|\n\psi|^2)
  \end{equation}
such that $ \rho > 1 $ if  $ |\n\psi|^2 \in [0, 1) $, therefore $  H  $ is a positive decreasing function defined on $[0, 1] $, twice
differentiable on $[0, 1) $, and satisfies $  H(1) = 1 $. 
Then \eqref{1-1}  reduces to a single equation
\begin{equation}\label{1-8}
\div \left(\frac{\n\psi}{H(|\n\psi|^2)}\right)=0.
\end{equation}

 \par   By using the hodograph method, Courant and Friedrichs in Section 104 of \cite{CF48} obtained some particular planar radially symmetric flows including circulatory flows and purely radial flows. Their superpositions are called
spiral flows. Weng-Xin-Yuan \cite{WXY20} gave a complete classification of the radially symmetric flow with or without shocks by prescribing suitable boundary conditions on  the inner and outer circle of an annulus and also analyzed the dependence of the solutions on the boundary data.  However, there are very few papers working on subsonic spiral flows outside a body. This motivate us to  investigate  the subsonic spiral flows whose asymptotic state is the radially symmetric subsonic spiral flow.  To this end, we first describe the subsonic spiral flow in the exterior domain $ B_1^c=\{(x_1,x_2):r=|x|>1\}$. It is convenient to use the polar coordinate $ (r,\theta) $:
\begin{equation}\label{1-9}
\begin{cases}
\begin{aligned}
&x_1=r\cos\theta,\\
&x_2=r\sin\theta,\\
\end{aligned}
\end{cases}
\end{equation}
and  decompose $ \textbf{u}=U_{1}\mathbf{e}_r+U_{2}\mathbf{e}_\theta $ with
\begin{equation}\label{1-10}
 \textbf{e}_r = (\cos\theta ,\sin\theta
)^T, \ \textbf{e}_\theta  = ( - \sin \theta ,\cos \theta)^T.
\end{equation}
Then the flow is described by radially symmetric smooth functions of the form $\mathbf{ u_0}(x)=U_{1}(r)\mathbf{e}_r+U_{2}(r)\mathbf{e}_\theta $, $ \rho(x)=\rho_b(r) $ and $ p(x)= p_b(r) $, which  solves the following system
 \begin{equation}\label{1-11}
\begin{aligned}
 \begin{cases}
 (\rho_b U_1)^{\prime}+\frac{1}{r} \rho_b U_{1}=0,&\quad  r>1,\\
\ U_{1}U_{1}^{\prime}
+\frac{1}{\rho_b} p_b^{\prime}-\frac{U_{2}^2}{r} =0,&\quad  r>1,\\
U_{1}U_{2}^{\prime}+\frac{U_{1}U_{2}}{r}
=0, &\quad  r>1
\end{cases}
\end{aligned}
\end{equation}
with the boundary condition
\begin{equation*}
\rho_b(1)=\rho_{0}>0, \ U_1(1)=\kappa_1>0,\ U_2(1)=\kappa_2.
\end{equation*}
Here $ \rho_{0} $, $ \kappa_1 $ and $ \kappa_2 $ are chosen to satisfy \eqref{1-4}, that is
\begin{equation*}
\frac{1}{2}(\kappa_1^2+\kappa_2^2)+\frac{\rho_0^{\gamma-1}}{\gamma-1}
=\frac{\gamma+1}{2(\gamma-1)}.
\end{equation*}
Hence it is easy to see that
\begin{equation}\label{1-12}
\rho_b U_1=\frac{\rho_0\kappa_1}{r},\quad
U_{2}=\frac{\kappa_2}{r}.
\end{equation}
Therefore
\begin{equation}\label{1-13}
\mathbf{u_0}(x)=\left(\frac{ \rho_b^{-1}\rho_0\kappa_1x_1-\kappa_2x_2}{r^2},
\frac{\rho_b^{-1}\rho_0\kappa_1x_2+\kappa_2x_1}{r^2}\right).
\end{equation}
Furthermore, by  \eqref{1-6}, the corresponding stream functions $ \psi_{0} $, $ \psi_{10} $ and $ \psi_{20} $ are
\begin{equation}\label{1-14}
\begin{aligned}
&\psi_{10}(\theta)=\rho_{0}\kappa_1\theta,\quad
\psi_{20}(r)=-\kappa_2\int_{1}^{r}\rho_b(s)\frac{1}{s}\mathrm{d}s,\\
&\psi_0=\psi_{10}+\psi_{20}=
\rho_0\kappa_1\theta-\kappa_2\int_{1}^{r}\rho_b(s)\frac{1}{s}\mathrm{d}s.
\end{aligned}
\end{equation}
\par Denote $ M_1^2=\frac{U_1^2}{c^2} $, $ M_2^2=\frac{U_2^2}{c^2}$ and $ M^2=M_1^2+M_2^2 $. Then by simple calculations, we have
\begin{equation}\label{1-15}
\begin{cases}
\begin{aligned}
&\frac{\de}{\de r}(M_1^2)
=-\frac{(2(1+M_2^2)+(\gamma-1)M^2)
M_1^2}{r(1-M_1^2)},\\
&\frac{\de}{\de r}(M_2^2)
=-\frac{(2(1-M_1^2)+(\gamma-1)M^2)
M_2^2}{r(1-M_1^2)},\\
&\frac{\de}{\de r}(M^2)
=-\frac{[(\gamma-1)M^2+2)]M^2 }{r(1-M_1^2)}.
\end{aligned}
\end{cases}
\end{equation}
Assume that $ \kappa_1^2+\kappa_2^2>1 $ and $ \kappa_1^2<1 $. Then  there exists  a smooth transonic spiral flows  for all $ r>1 $. The structural stability of this special transonic flows was investigated by Weng-Xin-Yuan in \cite{WYX20}.  Assume that
$ \kappa_1^2+\kappa_2^2<1 $. Then it follows from \eqref{1-15}  that $ M^2<1 $ and $ M_2^2<1 $ hold for all $ r>1 $, which means that the flow is uniformly subsonic.
\par Let $\md(\Gamma) $ is a porous body in $ \mathbb{R}^2 $, which contains the domain $ B_1(0) $. $ \Gamma $ is a bounded and connected  $ C^{\infty} $ smooth surface describing the boundary of $ \md $.  Suppose $ \Omega $ is  the exterior domain of $\md(\Gamma) $, $ i.e.,  \Omega:=\mathbb{R}^2 \backslash\md $,  which is connected and filled with compressible and invisible fluid.
\par In this paper, we aim to construct a smooth subsonic spiral flows outside  $\md(\Gamma) $ which tends to the radially symmetric subsonic spiral flow described above.  Since $\md(\Gamma) $ is a porous body, which means that on the boundary, we impose
\begin{equation}\label{1-16}
\rho \mathbf{u}\cdot \vec{n}=\rho_b(r) U_1(r)\mathbf{e}_r\cdot \vec{n},
\end{equation}
where $  \vec{n} $ stands for the unit inward normal of domain $ \md(\Gamma) $. It follows from \eqref{1-16} that $ \psi=\psi_{10} $   on $ \Gamma $.
\par Moreover, the radially symmetric subsonic spiral flow $ \mathbf{u_0}(x)$ decays like $ O(|x|^{-1}) $, then we may have the  decay rate  $|\mathbf{u}(x)- \mathbf{u_0}(x)|=O(|x|^{-2}) $ as $ |x|\rightarrow\infty $. In terms of the stream function,  we expect
\begin{equation}\label{1-17}
|\n (\psi-\psi_0)|=O(|x|^{-2}).
\end{equation}
\par Within this paper, we will consider the following problem:
\par $ \mathbf{Problem}$ $(\kappa_1,\kappa_2) $:  Find function $ \psi $ which  solving the following system:
\begin{equation}\label{1-18}
\begin{cases}
\div \left(\frac{\n\psi}{H(|\n\psi|^2)}\right)=0, &\quad \rm{in} \ \Omega,\\
\psi=\psi_{10}, &\quad \rm{on} \ \Gamma,\\
|\n (\psi-\psi_0)|=O(|x|^{-2}).\\
\end{cases}
\end{equation}
\par The main results of this paper can be stated as follows:
\begin{theorem}
 For any fixed $ \kappa_1\in(0,1) $, there exists a positive number $ \hat \kappa_2 $ such that for $ |\kappa_2| \in  [0,\hat \kappa_2)  $, there exists  a unique uniformly subsonic spiral flows to $ \textbf{Problem}$ $(\kappa_1,\kappa_2) $. More precisely, there exists a unique smooth solution $ \psi \in C^{\infty}(\Omega) $ to  \eqref{1-18} such that
 \begin{equation}\label{1-19}
\sup_{x\in\Omega}|\n \psi|<1.
\end{equation}
\end{theorem}
 \begin{theorem}
 Let $ \kappa_2^\epsilon\rightarrow \hat \kappa_2 $ as $ \epsilon\rightarrow 0 $ with  $ \kappa_2^\epsilon<\hat \kappa_2 $. Denote by $( u_1^\epsilon, u_2^\epsilon ) $ the uniformly subsonic spiral flow corresponding  to $ \textbf{Problem} $ $(\kappa_1,\kappa_2^\epsilon) $. Then there exists a subsequence, still labeled by $ ( u_1^\epsilon, u_2^\epsilon ) $ such that
\begin{equation}\label{1-20}
\quad u_1^\epsilon \rightarrow \hat{u}_1, \quad u_2^\epsilon \rightarrow \hat{u}_2,
\end{equation}
\begin{equation}\label{1-21}
g((q^\epsilon)^2) u_1^\epsilon  \rightarrow g(\hat q^2)\hat{u}_1, \quad
g((q^\epsilon)^2) u_2^\epsilon  \rightarrow g(\hat q^2)\hat{u}_2,
\end{equation}
where $ (q^\epsilon)^2= (u_1^\epsilon)^2+(u_2^\epsilon)^2 $, $ \hat q^2=\hat{u}_1^2+
\hat{u}_2^2 $, and $ g(q^2) $ is the function defined by \eqref{1-5}. All the above convergence are almost everywhere convergence. Moreover, this limit yields a subsonic-sonic spiral flow $ (\hat\rho,\hat u_1,\hat u_2) $, where $ \hat\rho=g(\hat q^2) $, which is a weak solution of $ \textbf{Problem} (\kappa_1,\hat\kappa_2) $.  The limit solution  $ (\hat\rho,\hat{u}_1, \hat{u}_2) $ satisfies \eqref{1-1}-\eqref{1-2} in the  sense of distribution and
 the boundary condition \eqref{1-16}  as the normal trace of the divergence-measure field $ (\hat\rho\hat u_1,\hat\rho\hat u_2) $ on the boundary.
\end{theorem}
\begin{theorem}
 For any fixed $ \kappa_2\in(-1,1) $, there exists a positive number $ \tilde \kappa_1 $ such that for  $ \kappa_1  \in  (0,\tilde \kappa_1)  $, there exists  a unique uniformly subsonic spiral flows to $ \textbf{Problem}$ $(\kappa_1,\kappa_2) $. More precisely, there exists a unique smooth solution $ \psi \in C^{\infty}(\Omega) $ to  \eqref{1-18} such that
  $\sup_{x\in\Omega}|\n \psi|<1 $.

 \end{theorem}
\begin{theorem}
 Let  $ \kappa_1^\epsilon\rightarrow \tilde \kappa_1 $ as $ \epsilon\rightarrow 0 $ with $ \kappa_1^\epsilon<\tilde \kappa_1 $.  Denote by $( u_1^\epsilon, u_2^\epsilon ) $ the uniformly subsonic spiral flow corresponding  to $ \textbf{Problem} $ $(\kappa_1^\epsilon,\kappa_2) $. Then there exists a subsequence, still labeled by $ ( u_1^\epsilon, u_2^\epsilon ) $ such that
\begin{equation}\label{1-22}
u_1^\epsilon \rightarrow \tilde{u}_1, \quad u_2^\epsilon \rightarrow \tilde{u}_2,
\end{equation}
\begin{equation}\label{1-23}
g((q^\epsilon)^2)\tilde{u}_1  \rightarrow g(\tilde q^2)\tilde{u}_1, \quad
g((q^\epsilon)^2)\tilde{u}_2  \rightarrow g(\tilde q^2)\tilde{u}_2,
\end{equation}
where $ (q^\epsilon)^2= (u_1^\epsilon)^2+(u_2^\epsilon)^2 $, $ \tilde q^2=\tilde{u}_1^2+
\tilde{u}_2^2 $.  All the above convergence are almost everywhere convergence. Moreover, this limit yields a subsonic-sonic spiral flow $ (\tilde\rho,\tilde u_1,\tilde u_2) $, where $ \tilde\rho=g(\tilde q^2) $, which is a weak solution of $ \textbf{Problem} (\tilde\kappa_1,\kappa_2) $.   The limit solution $ (\tilde\rho,\tilde{u}_1, \tilde{u}_2) $ satisfies \eqref{1-1}-\eqref{1-2} in the  sense of distribution and the boundary condition \eqref{1-16}  as the normal trace of the divergence-measure field $ (\tilde\rho\tilde u_1,\tilde\rho\tilde u_2) $ on the boundary.
\end{theorem}
\begin{remark}
The decay rate $ |\n (\psi-\psi_0)|=O(|x|^{-2}) $ is optimal as was observed by the results in \cite{CL11}, but here we remove the small perturbation conditions in \cite{CL11}. We employe the Kevin's transformation and some elliptic estimates to derive it.
\end{remark}

\begin{remark}
We prescribe the boundary condition \eqref{1-16} on the boundary $\Gamma$, not the one that $\rho \mathbf{u}\cdot \vec{n}=(\rho_b(r) U_1(r)\mathbf{e}_r+\rho_b(r) U_2(r){\bf e}_{\theta})\cdot \vec{n}$. Otherwise, the solution will be the background radially symmetric solution. One may pose $\rho \mathbf{u}\cdot \vec{n}=\nabla^{\perp} \Psi\cdot \vec{n}$ on $\Gamma$ for any given smooth function $\Psi$, our method still works in this case.
\end{remark}

\par The research on compressible inviscid flows has a long history, which provides many significant and challenging problems. The flow past a body, through a nozzle, and past a wall are typical flows patterns, which have physical significances and physical effects. The first theoretical result on the problem for irrotational flows past a body was obtained by Frankl and Keldysh in \cite{FK34}. The important progress for  two dimensional subsonic irrotational flows past a smooth body with a small free stream Mach number was obtained by Shiffman \cite{SM52}. Later on, Bers \cite{BL54} proved the existence of two dimensional subsonic irrotational flows around a profile with a sharp trailing edge and also showed that the maximum of Mach numbers approaches one as the free stream Mach number approaches the critical value. The uniqueness and asymptotic behavior of subsonic irrotational plane flows were studied by Finn and Gilbarg in \cite{RD57}.  The well-posedness theory for two-dimensional subsonic flows past a wall or a symmetric body was established by Chen, Du, Xie and Xin in \cite{CDXX16}. For the three-(or higher-) dimensional cases, the existence and uniqueness of three dimensional subsonic irrotational flows around a smooth body were proved by Finn and Gilbarg  in \cite{RF57}.  Dong and Ou \cite{DO93} extended the results of Shiffman to higher dimensions  by the direct method of calculus of variations and the standard Hilbert space method. The respective incompressible case is considered in \cite{OB94}. For the subsonic flow problem in nozzles, one may referee to \cite{DW18,DXY11,DXX14,GW17,LY14,WSK14, XC07,XC10}.
\par On the other hand, the subsonic-sonic limit solution can be constructed by the compactness method. The first  compactness framework on sonic-subsonic
irrotational flows in two dimension was introduced by  \cite{CDS07} and  \cite{XC07} independently. In \cite{CDS07}, Chen, Dafermos, Slemrod and Wang introduced general compactness framework and proved the existence of two-dimensional subsonic-sonic irrotational flows.   Xie and Xin  \cite{XC07}  established the subsonic-sonic limit of the two-dimensional irrotational flows
through infinitely long nozzles. Later on, they extended the result to the three-dimensional
axisymmetric flow through an axisymmetric nozzle in \cite{XC10}.    Furthermore, Huang, Wang and Wang \cite{HWW11} established a compactness framework for the general
multidimensional irrotational case.    Chen, Huang and Wang \cite{CHW16} established the compactness framework for nonhomentropic and rotational flows and proved the existence of multidimensional subsonic-sonic full Euler flows through infinitely long nozzles. Recently, by the compactness framework in \cite{CHW16}, the existence  of  subsonic-sonic flows with general conservative forces in an exterior domain was established by Gu and Wang in \cite{GX20}.
\par The rest of this article is organized as follows. In Section 2, we first introduce  the subsonic truncation to reformulate the problem into a second-order uniformly elliptic equation, and then establish the existence and uniqueness  of the modified  spiral flow  by a variational method. Finally we remove the truncation and complete the proof of Theorem 1.1 and Theorem 1.3. In Section 3,  the compensated compactness framework for steady irrotational flows is employed to establish the existence of weak subsonic-sonic spiral flows.
\section{Subsonic spiral flows }\noindent
\par This section is mainly devoted to the proof of Theorem 1.1 and Theorem 1.3. The proof can be divided into 5 subsections.
\subsection{Subsonic truncation}\noindent
\par  By direct calculations, it is easy to see that the derivative of function $ H(s) $ goes to  infinity as $ s\rightarrow 1 $. To control the ellipticity and avoid singularity of $ H^{\prime} $, one may truncate $ H $ as follows
\begin{equation}\label{2-1}
\tilde H(s)=
\begin{cases}
H(s), &\quad \rm{if}\ 0\leq s\leq 1-2\varepsilon,\\
\rm{smooth \ and \ decreasing}, &\quad \rm{if}\ 1-2\varepsilon\leq s\leq 1-\varepsilon ,\\
H(1-\varepsilon), &\quad \rm{if}\ s\geq 1-\varepsilon,
\end{cases}
\end{equation}
where $ \varepsilon $  is a small positive constant and $ \tilde H $ is a smooth decreasing function.

\par Then we consider the modified equation:
\begin{equation}\label{2-2}
\begin{cases}
\div \left(\frac{\n\psi}{\tilde H(|\n\psi|^2)}\right)=0, &\quad \rm{in} \ \Omega,\\
\psi=\psi_{10}, &\quad \rm{on} \ \Gamma,\\
|\n (\psi-\psi_0)|=O(|x|^{-2}).\\
\end{cases}
\end{equation}
After  the straightforward computation, $\eqref{2-2}_1 $ can be rewritten as
 \begin{equation}\label{2-3}
 \sum_{i,j=1}^{2}a_{ij}\p_{ij}\psi=0,
  \end{equation}
  where
  \begin{equation*}
a_{ij}=\frac{\tilde H\delta_{ij}-2\tilde H^{\prime}\p_i\psi\p_j\psi}{\tilde H^2}.
  \end{equation*}
  Then it is easy to verify  that there exist two constants $ \lambda $ and $ \Lambda $ such that
  \begin{equation}\label{2-4}
  \lambda|\xi|^2\leq a_{ij}\xi_i\xi_j\leq \Lambda|\xi|^2, \quad {\rm{for}}\ \xi\in \mathbb{R}^2.
\end{equation}
Here $ \lambda $ and $ \Lambda $  depend on $ \gamma $ and $ \varepsilon $.
\par Next, we solve \eqref{2-2} by a variational method.  We first give the definition of weak solution to be used  in next subsection.
\begin{definition}
A function $ \psi\in H_{loc}^1(\Omega)  $ is a weak solution to \eqref{2-2} if
\begin{equation*}
\int_{\Omega}(\tilde H(|\nabla\psi|^2))^{-1}\nabla\psi\cdot \nabla v\mathrm{d}x=0
\end{equation*}
holds for any $ v\in C_c^{\infty}(\Omega) $.
\end{definition}

\subsection{Existence of weak solution of the problem (2.2)}\noindent

\par Define the space
\begin{equation}\label{2-5}
\ml=\{\phi\in  L_{loc}^2(\Omega),\nabla\phi\in  L^2(\Omega),\phi|_{\Gamma}=0\}.
\end{equation}
It is easy to see that $ \ml $ is a Hilbert space under the $ \mathring{H}^1 $ norm. One should look for the solution $ \psi $ to the problem \eqref{2-2} with the form: $ \psi=\phi_0+\psi_0 $, where $ \phi_0\in L_{loc}^2(\Omega) $ and $ \nabla\phi_0\in  L^2(\Omega) $. It follows from the boundary condition in \eqref{2-2} that  $ \phi_0=-\psi_{20} $  on $ \Gamma $. To homogenize the boundary data on $ \Gamma $, we should introduce another function $ \varphi_0 $ such that $ \phi\triangleq\phi_0-\varphi_0 \in \ml  $.
\par Let $ \zeta(x) $ be a cut-off function such that $ \zeta\equiv1 $ near $ \Gamma $, $ \zeta\equiv 0 $ outside a circle $ B_{R_0} $ with $ R_0 $ large enough. Then set  $ \varphi_0 =-\zeta\psi_{20} $. By simple calculations,  $ -\zeta\psi_{20}$  satisfies all the requirements we imposed. 
\par Define
\begin{equation*}
\psi=\phi+\varphi_0+\psi_0, \quad
F(s)=\int_{0}^{s}(\tilde H(t))^{-1}\mathrm{d}t,
\end{equation*}
and
\begin{equation}\label{2-6}
I[\phi;\varphi_0,\psi_0]=\int_{\Omega}
[F(|\nabla{\phi}+\nabla\varphi_0+\nabla\psi_0|^2)
-F(|\nabla\psi_0|^2)
-2F^{\prime}(|\nabla\psi_0|^2)\nabla\psi_0\cdot (\nabla{\phi}+\nabla\varphi_0)]\mathrm{d}x.\\
\end{equation}
We consider the following variational problem:
\begin{equation}\label{2-7}
l=\min_{\phi\in{\ml}}I[\phi;\varphi_0,\psi_0].
\end{equation}

 \par For our variational problem, we have the following theorem:
 \begin{theorem}
 The functional $ I $ has  uniquely one minimizer $ \phi  $, i.e. $ I[\phi;\varphi_0,\psi_0]=l $, then $ \psi=\phi+\varphi_0+\psi_0 $ satisfies
 \begin{equation}\label{2-8}
\begin{cases}
\div \left(\frac{\n\psi}{\tilde H(|\n\psi|^2)}\right)=0, &\quad \rm{in} \ \Omega,\\
\psi=\psi_{10}, &\quad \rm{on} \ \Gamma,\\
\end{cases}
\end{equation}
in the weak sense.  In addition, there is a constant $ C $ such that
\begin{equation}\label{2-9}
\int_{\Omega}|\n\phi|^2\mathrm{d}x\leq C.
\end{equation}
The constant  $ C $ here and in the rest of this paper depends on $ \Omega $, $ \gamma $, $\varepsilon $, $ \kappa_1 $ and $ \kappa_2$. We will not repeat the dependence every time.
\end{theorem}
$\mathbf{Proof} $. $\textbf{Step 1} $. $ I[\phi;\varphi_0,\psi_0] $ is coercive on $ \ml $. Let $ p=(p_1,p_2) $, $ \mf(p)=F(|p|^2) $.  It follows from \eqref{2-4} that
\begin{equation}\label{2-10}
 \forall \xi\in \mathbb{R}^2,\ \lambda|\xi|^2\leq \p_{p_ip_j}^2\mf(p)\xi_i\xi_j\leq \Lambda|\xi|^2.
\end{equation}
Then by direct calculation, we have
\begin{equation*}
\begin{aligned}
&F(|\n\phi+\n\varphi_0+\n\psi_0|^2)-F(|\n\psi_0|^2)
-2F^{\prime}(|\n\psi_0|^2)\n\psi_0\cdot (\n\phi+\n\varphi_0)\\
&=\mf(\n\psi)-\mf(\n\psi_0)-\p_p\mf(\n\psi_0)\cdot (\n\phi+\n\varphi_0)\\
&=\int_{0}^{1}t\p_{p_ip_j}^2\mf(\n\psi_0+(1-t)(\n\phi+\n\varphi_0))
\p_i(\phi+\varphi_0) \p_j(\phi+\varphi_0)\mathrm{d}t.\\
\end{aligned}
\end{equation*}
This implies
\begin{equation}\label{2-11}
\Lambda\|\phi+\varphi_0\|_{\ml}^2\geq I \geq \lambda\|\phi+\varphi_0\|_{\ml}^2\geq \frac{\lambda}{2}\|\phi\|_{\ml}^2
-\lambda\int_{\Omega}|\n\varphi_0|^2\mathrm{d}x.
\end{equation}
Therefore we obtain
\begin{equation}\label{2-12}
I\geq \frac{\lambda}{2}\|\phi\|_{\ml}^2-C.
\end{equation}
From \eqref{2-12}, one can conclude that the energy functional is coercive on $ \ml $.

\par $ \textbf{Step 2} $. The existence of the minimizer $ \phi \in \ml $. By coerciveness of the functional $ I[\cdot;\varphi_0,\psi_0] $ on $ \ml $, we know that $ l=\inf_{\phi\in\ml}I[\phi;\varphi_0,\psi_0] $ exists. Since $ I[0;\varphi_0,\psi_0]\leq \Lambda\|\varphi_0\|_{\ml}^2\leq C $, we have $ l\leq C $. By  definition, there is a sequence $ \{\phi_k\}_{k=1}^{\infty} \subset \ml $ such that
$ \lim_{k\rightarrow \infty}I[\phi_k;\varphi_0,\psi_0]=l $. For sufficiently large $ k $, $ C\geq I[\phi_k;\varphi_0,\psi_0]\geq \frac{\lambda}{2}\|\phi_k\|_{\ml}^2-C $,  hence $ \|\phi_k\|_{\ml}^2\leq \frac{4C}{\lambda} $. Since $ \ml $ is a Hilbert space, there exists a subsequence of $ \{\phi_k\} $ (still denoted by $ \{\phi_k\} $), that converges weakly to some $ \phi \in \ml $.
\par Next, we need to show that the functional $ I[\cdot;\varphi_0,\psi_0] $ is lower semi-continuous with respect to the weak convergence of $ \ml $, that is
\begin{equation}\label{2-13}
I[\phi;\varphi_0,\psi_0]\leq \liminf_{k\rightarrow\infty}
I[\phi_k;\varphi_0,\psi_0],\ {\rm{if}}\ \phi_k\rightharpoonup \phi\ \rm{in}\ \ml.
\end{equation}
  It follows from  $\phi_k\rightharpoonup \phi\ \rm{in}\ \ml $ that $ \n\phi_k\rightharpoonup \n\phi  $ in $  L_{loc}^{2}(\Omega) $. Furthermore, one has  $ 2F^{\prime}(|\n\psi_0|^2)\n\psi_0 \in  L_{loc}^{2}(\Omega) $. Denote $ \Omega^{R_j}=\Omega\bigcap B_{R_j}(0)$, where $ \{R_j\}_{j=1}^{\infty}\in \mathbb{R} $ is  an increasing sequence  and $ R_j\rightarrow \infty $ as $ j\rightarrow \infty $. Then
\begin{equation}\label{2-14}
\lim_{k\rightarrow \infty}\int_{\Omega^{R_j}}2F^{\prime}(|\n\psi_0|^2)\n\psi_0
\cdot(\n\phi_k+\n\varphi_0)\mathrm{d}x
=\int_{\Omega^{R_j}}2F^{\prime}(|\n\psi_0|^2)\n\psi_0
\cdot(\n\phi+\n\varphi_0)
\mathrm{d}x.
\end{equation}
For $ \psi_k=\phi_k+\varphi_0+\psi_0 $, a simple calculation yields that
\begin{equation}\label{2-15}
\begin{aligned}
&F(|\n\psi_k|^2)-F(|\n\psi|^2)
-2F^{\prime}(|\n\psi|^2)\n\psi\cdot (\n\phi_k-\n\phi)\\
&=\int_{0}^{1}t\p_{p_ip_j}^2\mf(\n\psi+(1-t)\n(\phi_k-\phi))\p_i(\phi_k-\phi) \p_j(\phi_k-\phi)\mathrm{d}t\geq 0.
\end{aligned}
\end{equation}
Similar to \eqref{2-14}, we have
\begin{equation*}
\lim_{k\rightarrow \infty}  \int_{\Omega^{R_j}}2F^{\prime}(|\n\psi|^2)\n\psi\cdot (\n\phi_k-\n\phi)=0.
\end{equation*}
So integrating both sides of \eqref{2-15} over $ \Omega^{R_j} $ leads to
\begin{equation}\label{2-16}
\liminf_{k\rightarrow \infty}  \int_{\Omega^{R_j}}F(|\n\psi_k|^2)-F(|\n\psi|^2)\mathrm{d}x\geq 0.
\end{equation}
Collecting \eqref{2-14}-\eqref{2-15} and \eqref{2-16} together  gives
\begin{equation}\label{2-17}
\begin{aligned}
 &\liminf_{k\rightarrow \infty} \int_{\Omega^{R_j}}[F(|\n\psi_k|^2)-F(|\n\psi_0|^2)
-2F^{\prime}(|\n\psi_0|^2)\n\psi_0\cdot (\n\phi_k+\n\varphi_0)]\mathrm{d}x\\
&\quad \quad \geq \int _{\Omega^{R_j}}[F|\n\psi|^2)-F(|\n\psi_0|^2)
-2F^{\prime}(|\n\psi_0|^2)\n\psi_0\cdot (\n\phi+\n\varphi_0)]\mathrm{d}x.
\end{aligned}
\end{equation}
It follows from \eqref{2-17} that we deduce
\begin{equation}\label{2-18}
\liminf_{k\rightarrow\infty}I[\phi_k,\varphi_0,\psi_0]
\geq \int_{\Omega^{R_j}}F(|\n\psi|^2)-F(|\n\psi_0|^2)
-2F^{\prime}(|\n\psi_0|^2)\n\psi_0\cdot (\n\phi+\n\varphi_0)]\mathrm{d}x.
\end{equation}
This inequality holds for each $ R_j $. Let $ R_j\rightarrow \infty $ and use the monotone convergence theorem to conclude
\begin{equation*}
\liminf_{k\rightarrow\infty}I[\phi_k;\varphi_0,\psi_0]\geq I[\phi;\varphi_0,\psi_0].
\end{equation*}
Therefore we proved \eqref{2-13}. $ \phi \in \ml $ is actually a minimizer.
 \par $\textbf{Step 3} $. The uniqueness of the minimizer $ \phi $.  For any $ \phi_1,\phi_2 \in \ml $, we derive that
\begin{equation}\label{2-19}
\begin{aligned}
&I[\phi_1;\varphi_0,\psi_0]
+I\left[\phi_2;\varphi_0,\psi_0\right]
-2I[\frac{\phi_1+\phi_2}{2};\varphi_0,\psi_0]\\
&=\int_{\Omega}\mf(\n\psi_1)-\mf(\n\psi_2)
-2\mf\left(\frac{\n\psi_1+\n\psi_2}{2}\right)\mathrm{d}x\\
&\geq\frac{\lambda}{2}\|\psi_1-\psi_2\|_{\ml}^2
=\frac{\lambda}{2}\|\phi_1-\phi_2\|_{\ml}^2.
\end{aligned}
\end{equation}
For the uniqueness, suppose $ \phi_1 $ and $ \phi_2 $ are two minimizers. Then $\frac{\phi_1+\phi_2}{2}\in\ml $, and
\begin{equation}\label{2-20}
I[\phi_1;\varphi_0,\psi_0]+I[\phi_2;\varphi_0,\psi_0]\leq 2I\left[\frac{\phi_1+\phi_2}{2};\varphi_0,\psi_0\right].
\end{equation}
 Using \eqref{2-19}, one has
\begin{equation}\label{2-21}
0\geq I[\phi_1;\varphi_0,\psi_0]+I[\phi_2;\varphi_0,\psi_0]- 2I\left[\frac{\phi_1+\phi_2}{2};\varphi_0,\psi_0\right]\geq
\frac{\lambda}{2}\|\phi_1-\phi_2\|_{\ml}^2,
\end{equation}
 which implies that the minimizer is unique in $ \ml $.
\par $\textbf{Step 4} $. $ \psi=\phi+\varphi_0+\psi_0 $ satisfies \eqref{2-8} in the weak sense.  It follows from  $ \textbf{Step 1-3} $ that the functional $ I $ has a unique minimizer  $ \phi $.  For any  $ v\in C_c^{\infty}(\Omega) $, define $ I[\phi+\tau v;\varphi_0,\psi_0] $ for $ \tau>0 $. Then
\begin{equation*}
\begin{aligned}
I(\tau)-I(0)&=\int_{\Omega}[F(\n\phi +\tau \n v+\n\varphi_0+\n\psi_0|^2)-F(|\n\phi+\n\varphi_0
+\n\psi_0|^2)]\mathrm{d}x\\
 &\quad -\tau\left(\int_{\Omega}2F^{\prime}(|\n\psi_0|^2)\n\psi_0\cdot \n v\mathrm{d}x\right)\\
&=V_1+\tau V_2.
\end{aligned}
\end{equation*}
Since $ \div (F^{\prime}(|\n\psi_0|^2)\n\psi_0)=0 $, so $ V_2=0 $ for any $ v\in C_c^{\infty}(\Omega) $. Hence we deduce that
\begin{equation*}
\begin{aligned}
\frac{I(\tau)-I(0)}{\tau}&=\frac{1}{\tau}\int_{\Omega}[F(|\n\phi +\tau \n v+\n\varphi_0+\n\psi_0|^2)
-F(|\n\phi+\n\varphi_0+\n\psi_0|^2)]\mathrm{d}x\\
&=\int_{\Omega}\int_{0}^{1}F^{\prime}(|\n\psi +s\tau \n v|^2)(\n \psi +s\tau \n v)\cdot \n v
\mathrm{d}s\mathrm{d}x\\
&=\int_{\Omega}\int_{0}^{1}(F^{\prime}(|\n\psi +s\tau \n v|^2)-F^{\prime}(|\n \psi|^2))\n\psi\cdot \n v
\mathrm{d}s\mathrm{d}x\\
 &\quad+\tau\int_{\Omega}\int_{0}^{1}F^{\prime}(|\n \psi +s\tau \n v|^2)s|\n v|^2\mathrm{d}s\mathrm{d}x+\int_{\Omega}F^{\prime}(|\n\psi|^2)\n\psi\cdot \n v\mathrm{d}s\mathrm{d}x\\
&=W_1+W_2+W_3.
\end{aligned}
\end{equation*}
 Due to truncation, $ F^{\prime} $ is bounded. Then $ W_2\leq\tau \Lambda\|v\|_{\ml}^2 $, so $ \lim_{\tau\rightarrow 0}W_2=0 $. By Lebesgue convergence theorem, we get $ \lim_{\tau\rightarrow 0}W_1=0 $. Since $ I^{\prime}(0)=0 $, this gives $ W_3=0 $. Therefore $ \psi=\phi+\varphi_0+\psi_0 $ is a weak solution to \eqref{2-8}.
\subsection{Regularity of weak solution}\noindent
\par In this section, we verify that the weak solution obtained from the variation in subsection 2.2 is classical solution.
\begin{lemma}
  Suppose that $ \psi\in H_{loc}^1(\Omega)  $ is a weak solution to problem \eqref{2-2}. Then there is a constant $ \alpha\in(0,1) $ such that for any subregion $ \Omega_1\Subset (\Gamma\cup\Omega)$, there holds
\begin{equation}\label{2-22}
\sup_{x\in \Omega_1 }|\n\psi|+
\sup_{x_1,x_2\in \Omega_1}\frac{|\n\psi(x_1)-\n\psi(x_2)|}{|x_1-x_2|^{\alpha}}\leq C,
\end{equation}
where $C$ depends on $\Omega_1,\Omega,\gamma, \epsilon,\rho_0,\kappa_1,\kappa_2$.
\end{lemma}
$ \mathbf{Proof} $: Denote $ \bar\psi_k=\p_k \psi $ for $ k=1,2 $. It is easy to verify that
 \begin{equation}\label{2-23}
 \p_i(a_{ij}\p_j\bar\psi_k)=0,
 \end{equation}
 where $ a_{ij}=F^{\prime}(|\n\psi|^2)\delta_{ij}+2F^{''}(|\n\psi|^2)\p_i\psi \p_j\psi $. $ [a_{ij}] $ has uniformly positive eigenvalues.
 \par For the interior estimate, let $ U $  be a bounded interior subregion of $ \Omega $. It follows from  Theorem 8.24 in \cite{DN83} that $ \|\bar\psi_k\|_{C^{\alpha}(U)}\leq C\|\bar\psi_k\|_{L^2(\Omega)} $. According to \eqref{2-9}, we know that $ \|\bar\psi_k\|_{L^2(\Omega)} $ is uniformly bounded. From the standard elliptic estimate, we obtain \eqref{2-22}.
 \par Next for the boundary estimate near $ \Gamma $, one can apply Theorem 8.29 in \cite{DN83} and follow the arguments above to show \eqref{2-22} holds.
 \par Then a standard bootstrap argument implies that $ \psi $ actually belongs to $ C^{\infty}(\Omega) $ and $ \psi $ is a classical solution.

  \par Now, we use $ L^{\infty} $ bound of gradient to show that $ |\nabla(\psi-\psi_0)| $
  tend to zero as $ |x|\rightarrow\infty $ with decay rate $ |x|^{-2}$ in $ \Omega $.
  \begin{lemma}
 There is an  estimate of the continuity of $ \n\psi $ at infinity:
 \begin{equation}\label{2-24}
 |\n\psi-\n\psi_0|\leq \frac{ C}{(1+|x|)^{2}}.
 \end{equation}
 \end{lemma}
  $ \mathbf{Proof} $: For a fixed $ R\gg 1 $ with $ N=\{|x|\geq R\}\subset\Omega $.
   Set $ \bar\psi_1=\p_1 \psi=-\rho u_2 $ and $\bar\psi_2=\p_2 \psi=\rho u_1 $. Then we have the following elliptic system
  \begin{equation}\label{2-25}
  -a_{22}(x)\p_{x_2}\bar\psi_{2}=a_{11}(x)\p_{x_1}\bar\psi_{1}
  +2a_{12}(x)\p_{x_2}\bar\psi_{1}, \ \p_{x_1}\bar\psi_{2}= \p_{x_2}\bar\psi_{1},
 \end{equation}
 where
 \begin{equation*}
\begin{cases}
a_{22}(x)=F^{\prime}(|\n\psi|^2)+2F^{''}(|\n\psi|^2)\bar\psi_2^2,\\
a_{12}(x)=2F^{\prime}(|\n\psi|^2)+2F^{''}(|\n\psi|^2)\bar\psi_1\bar\psi_2,\\
a_{11}(x)=F^{\prime}(|\n\psi|^2)+2F^{''}(|\n\psi|^2)\bar\psi_1^2.\\
\end{cases}
\end{equation*}
  Applying Theorem 3 in \cite{RD57}, we know that $ \bar\psi_1 $ and $ \bar \psi_2 $ can be represented in the form
\begin{equation}\label{2-26}
\left(\begin{matrix}
\bar\psi_{1}\\
\bar\psi_{2}
\end{matrix}\right)
=\left(\begin{matrix}
\bar\psi_{10}\\
\bar\psi_{20}
\end{matrix}\right)+\frac{1}{c_0x_1^2+b_0x_1x_2+a_0x_2^2 }
\left(\begin{matrix}
\beta_1x_1+\beta_2x_2\\
\gamma_1x_1+\gamma_2x_2
\end{matrix}\right)
+\left(\begin{matrix}
O(|x|^{-1-\alpha})\\
O(|x|^{-1-\alpha})
\end{matrix}\right), \ {\rm{as}}\ |x|\rightarrow\infty,
\end{equation}
where $ 0<\alpha<1 $ and
\begin{equation*}
\begin{cases}
\begin{aligned}
&\lim_{|x|\rightarrow\infty}\bar\psi_1=\bar\psi_{10}, \quad \lim_{|x|\rightarrow\infty}\bar\psi_2=\bar\psi_{20},\\
&c_0=F^{\prime}((\bar\psi_{10}+\bar\psi_{20})^2)
+2F^{''}((\bar\psi_{10}+\bar\psi_{20})^2)\bar\psi_{20}^2,\\
&b_0=2F^{''}((\bar\psi_{10}+\bar\psi_{20})^2)\bar\psi_{10}\bar\psi_{20},\\
&a_0=F^{\prime}((\bar\psi_{10}+\bar\psi_{20})^2)
+2F^{''}((\bar\psi_{10}+\bar\psi_{20})^2)\bar\psi_{10}^2,\\
&c_0(\gamma_1+\beta_2)+2b_0\beta_1=0, \quad  c_0\gamma_2=a_0\beta_1.\\
\end{aligned}
\end{cases}
\end{equation*}
Since $ \n \psi-\n \psi_0 \in L^2(\Omega) $, then it follows from \eqref{2-26} that $ \bar\psi_{10}=\bar\psi_{20}=0 $. Thus $ a_0=c_0=F^{\prime}(0) $, $ b_0=0 $, $ \gamma_1+\beta_2=0 $ and $ \gamma_2=\beta_1 $. Hence \eqref{2-26} can be rewritten as
\begin{equation}\label{2-27}
\n \psi=\left(\begin{matrix}
\bar\psi_{1}\\
\bar\psi_{2}
\end{matrix}\right)
=\frac{1}{F^{\prime}(0)(x_1^2+x_2^2)}\left\{\beta_1
\left(\begin{matrix}
x_1\\
x_2
\end{matrix}\right)
+\beta_2
\left(\begin{matrix}
x_2\\
-x_1
\end{matrix}\right)\right\}+\left(\begin{matrix}
O(|x|^{-1-\alpha})\\
O(|x|^{-1-\alpha})
\end{matrix}\right), \ {\rm{as}}\ |x|\rightarrow\infty.
\end{equation}
Then $ \n \psi-\n \psi_0 \in L^2(\Omega) $ implies that
\begin{equation}\label{2-28}
\n \psi_0=\frac{1}{F^{\prime}(0)(x_1^2+x_2^2)}\left\{\beta_1
\left(\begin{matrix}
x_1\\
x_2
\end{matrix}\right)
+\beta_2
\left(\begin{matrix}
x_2\\
-x_1
\end{matrix}\right)\right\},\ {\rm{as}}\ |x|\rightarrow\infty.
\end{equation}
Since $ \psi=\phi_0+\psi_0 $, thus
 \begin{equation}\label{2-29}
\n\phi_0=\left(\begin{matrix}
O(|x|^{-1-\alpha})\\
O(|x|^{-1-\alpha})
\end{matrix}\right), \ {\rm{as}}\ |x|\rightarrow\infty.
\end{equation}
Furthermore, it follows from \eqref{2-25} that we derive
  \begin{equation}\label{2-30}
  \sum_{i,j=1}^{2}a_{ij}(x)\p_{x_ix_j}\bar\psi_{1}
  +\sum_{i=1}^{2}c_{i}(x)\p_{x_i}\bar\psi_{1}=0,
  \end{equation}
  where
  \begin{equation*}
   \begin{aligned}
  c_1(x)&=\p_{x_1}a_{11}(x)-\frac{a_{11}(x)}{a_{22}(x)}\p_{x_1}a_{22}(x),\\
   c_2(x)&=2\p_{x_1}a_{12}(x)-\frac{2a_{12}(x)}{a_{22}}\p_{x_1}a_{22}(x).
    \end{aligned}
\end{equation*}
 We set $ \bar\phi_{0}=\p_1\phi_0 $ and $  \bar\psi_{0}=\p_1 \psi_0 $. Then $ \bar\phi_{0} $  satisfies
\begin{equation}\label{2-31}
  \sum_{i,j=1}^{2}a_{ij}(x)\p_{x_ix_j}\bar\phi_{0}
  +\sum_{i=1}^{2}c_{i}(x)\p_{x_i}\bar\phi_{0}
   +\sum_{i,j=1}^{2}a_{ij}(x)\p_{x_ix_j}\bar\psi_{0}
  +\sum_{i=1}^{2}c_{i}(x)\p_{x_i}\bar\psi_{0}=0
  \end{equation}
  Introduce the Kelvin's transform
   \begin{equation}\label{2-32}
y_1=\frac{R^2}{|x|^2}x_1=:\Phi^1(x),\quad y_2=\frac{R^2}{|x|^2}x_2=:\Phi^2(x),
\end{equation}
 and write
 \begin{equation*}
y=\mathbf{\Phi}(x), \quad x=\mathbf{\Psi}(y),
 \end{equation*}
  where $ \mathbf{\Psi}=\mathbf{\Phi}^{-1} $.
  \par   Set  $  \phi_0^\sharp(y)=\bar\phi_{0}(\mathbf{\Psi}(y))$ and  $  \psi_0^\sharp(y)=\bar \psi_{0}(\mathbf{\Psi}(y))$. Then it follows from \eqref{2-29} and \eqref{2-32} that $ \phi_0^\sharp(y)=O(|y|^{1+\alpha}) $. Thus
  \begin{equation}\label{2-33}
   \phi_0^\sharp(0)=0,\quad
   \n\phi_0^\sharp(0)=0.
 \end{equation}
  Meanwhile, under this transformation, \eqref{2-31} can be written as
  \begin{equation}\label{2-34}
 \sum_{i,j=1}^{2}{a}_{ij}^\sharp(y) \p_{y_iy_j}\phi_{0}^\sharp+\sum_{i=1}^{2} {c}_{j}^\sharp(y) \p_{y_j}\phi_{0}^\sharp= f^\sharp(y),
 \end{equation}
 where
 \begin{equation*}
\begin{aligned}
   {a}_{11}^\sharp(y)&=a_{11}(\mathbf{\Psi}(y))(\Phi^1_{x_1})^2
   +2a_{12}(\mathbf{\Psi}(y))\Phi^1_{x_1}\Phi^1_{x_2}
   +a_{22}(\mathbf{\Psi}(y))(\Phi^1_{x_2})^2,\\
   {a}_{12}^\sharp(y)&= {a}_{21}^\sharp(y)=a_{11}(\mathbf{\Psi}(y))\Phi^1_{x_1}\Phi^2_{x_1}
   +a_{12}(\mathbf{\Psi}(y))(\Phi^1_{x_1}\Phi^2_{x_2}+\Phi^1_{x_2}\Phi^2_{x_1})
   +a_{22}(\mathbf{\Psi}(y))\Phi^1_{x_2}\Phi^2_{x_2},\\
  {a}_{22}^\sharp(y)&=a_{11}(\mathbf{\Psi}(y))(\Phi^2_{x_1})^2
   +2a_{12}(\mathbf{\Psi}(y))\Phi^2_{x_1}\Phi^2_{x_2}
   +a_{22}(\mathbf{\Psi}(y))(\Phi^2_{x_2})^2,\\
    {c}_1^\sharp(y)&=a_{11}(\mathbf{\Psi}(y))\Phi^1_{x_1x_1}+
   2a_{12}(\mathbf{\Psi}(y))\Phi^1_{x_1x_2}
   +a_{22}(\mathbf{\Psi}(y))\Phi^1_{x_2x_2}\\
   & \quad+\sum_{i=1}^{2}\left(\p_{y_i}a_{11}(\mathbf{\Psi}(y))
   -\frac{a_{11}(\mathbf{\Psi}(y))}{a_{22}(\mathbf{\Psi}(y))}
   \p_{y_i}a_{22}(\mathbf{\Psi}(y))\right)\Phi^1_{x_1}\Phi^i_{x_1}\\
  &\quad+2\sum_{i=1}^{2}\left(\p_{y_i}a_{12}(\mathbf{\Psi}(y))
   -\frac{a_{11}(\mathbf{\Psi}(y))}{a_{22}(\mathbf{\Psi}(y))}
   \p_{y_i}a_{12}(\mathbf{\Psi}(y))\right)\Phi^i_{x_1}\Phi^1_{x_2},\\
    {c}_2^\sharp(y)&=a_{11}(\mathbf{\Psi}(y))\Phi^2_{x_1x_1}+
   2a_{12}(\mathbf{\Psi}(y))\Phi^2_{x_1x_2}
   +a_{22}(\mathbf{\Psi}(y))\Phi^2_{x_2x_2}\\
   & \quad+\sum_{i=1}^{2}\left(\p_{y_i}a_{11}(\mathbf{\Psi}(y))
   -\frac{a_{11}(\mathbf{\Psi}(y))}{a_{22}(\mathbf{\Psi}(y))}
   \p_{y_i}a_{22}(\mathbf{\Psi}(y))\right)\Phi^2_{x_1}\Phi^i_{x_1}\\
  &\quad+2\sum_{i=1}^{2}\left(\p_{y_i}a_{12}(\mathbf{\Psi}(y))
   -\frac{a_{11}(\mathbf{\Psi}(y))}{a_{22}(\mathbf{\Psi}(y))}
   \p_{y_i}a_{12}(\mathbf{\Psi}(y))\right)\Phi^i_{x_1}\Phi^2_{x_2},\\
     f^\sharp(y)&= -\sum_{i,j=1}^{2} {a}_{ij}^\sharp(y) \p_{y_iy_j} \psi_{0}^\sharp-\sum_{i=1}^{2} {c}_{j}^\sharp(y)\p_{y_j} \psi_{0}^\sharp.
 \end{aligned}
\end{equation*}
Applying Schauder estimates in \cite{DN83}, we deduce that
 \begin{equation}\label{2-35}
|\phi_0^\sharp|_{2,\alpha;B_{{R}/{2}}}\leq  C(|\phi_0^\sharp|_{0;B_R}+| f^\sharp|_{0,\alpha;B_R}).
 \end{equation}
 Since $ \p_i(a_{ij}\p_j\bar\phi_0)+\p_i(a_{ij}\p_j\bar\psi_0)=0 $,  we can apply Theorem 8.15 in \cite{DN83} to derive that $ |\bar\phi_0|_{0;N}\leq C $.  By substituting $ |\phi_0^\sharp|_{0;B_R}=|\bar\phi_0|_{0;N}\leq C $ into \eqref{2-35}, we infer that
\begin{equation}\label{2-36}
|\phi_0^\sharp|_{2,\alpha;B_{R/2}}\leq  C.
\end{equation}
Combing \eqref{2-33}, we obtain
\begin{equation}\label{2-37}
|\phi_0^\sharp|\leq  C|y|^{2}, \quad  {\rm{in}}\ B_{R/2}.
\end{equation}
Hence
\begin{equation}\label{2-38}
|\bar \phi_0 | \leq \frac{ C}{|x|^{2}}, \quad  {\rm{in}}\ \{|x|\geq 2R\}.
\end{equation}
 Thus \eqref{2-24} can be  obtained from \eqref{2-38}.
\subsection{Modified spiral flows}\noindent
\par For the modified equation \eqref{2-2}, we have the following theorem:
\begin{theorem}
For every $ \psi_0 $,  \eqref{2-2} has a unique classical spiral  solution $ \psi $ such that $ \psi=\phi+\varphi_0+\psi_0 $ with $ \phi \in \ml $. Furthermore, for each fixed $ \kappa_1 $,  $\n \psi $ depends on $ \kappa_2 $ continuously. In particular, $ \max_{\Omega}|\n \psi| $ is a continuous function of $ |\kappa_2 | $. For each fixed $ \kappa_2 $, we have the similar conclusion.
\end{theorem}
$ \mathbf{Proof} $:  The existence follows from Theorem 2.2 and the regularity estimate follows  Lemma 2.3 and 2.4. To prove the uniqueness, we assume that  classical solutions
\begin{equation*}
\psi_i=\phi_i+\varphi_0+\psi_0,\quad i=1,2 \quad {\rm{with}}\ \phi_i\in\ml
\end{equation*}
would be both critical points of $ I[\phi;\varphi_0,\psi_0] $. Let $ I_{\phi}^{\prime} $ be the Fr$\acute{\rm{e}}$chet derivative.  Then we have
\begin{equation*}
\begin{aligned}
0&=(I_{\phi}^{\prime}(\phi_1;\varphi_0,\psi_0)
-I_{\phi}^{\prime}(\phi_2;\varphi_0,\psi_0),\phi_1-\phi_2)\\
&=\int_{\Omega}[\p_p\mf(\n\phi_1+\varphi_0+\psi_0)-
\p_p\mf(\n\phi_2+\n\varphi_0+\n\psi_0)](\n\phi_1-\n\phi_2)\mathrm{d}x\\
&=\int_{\Omega}\int_{0}^{1}\p_{p_ip_j}^2\mf(t\n\phi_1+(1-t)\n\phi_2
+\n\varphi_0+\n\psi_0)\p_i(\phi_1-\phi_2)\p_j(\phi_1-\phi_2)
\mathrm{d}t\mathrm{d}x\\
&\geq \lambda\|\phi_1-\phi_2\|_{\ml}^2.
\end{aligned}
\end{equation*}
Hence $ \phi_1=\phi_2 $.
\par  Now we examine the continuity of $ I[\phi; \varphi_0,\psi_0] $. Let $ \psi_i=\phi_i+\varphi_0+\psi_0 $, $ i=1,2 $, one has
\begin{equation*}
\begin{aligned}
&\mf(\n\phi_1+\n\varphi_0+\n\psi_0)-\mf(\n\psi_0)
-\p_p\mf(\n\psi_0)(\n\phi_1+\n\varphi_0)
-[\mf(\n\phi_2+\n\varphi_0+\n\psi_0)\\
&-\mf(\n\psi_0)-\p_p\mf(\n\psi_0)(\n\phi_2+\n\varphi_0)]\\
&=\int_{0}^{1}\p_p\mf(t\n\phi_1+(1-t)\n\phi_2+\n\varphi_0+\n\psi_0)\mathrm{d}t
(\n\phi_1-\n\phi_2)
-\p_p\mf(\n\psi_0)(\n\phi_1-\n\phi_2)\\
&=\int_{0}^{1}\int_{0}^{1}\p_{p_ip_j}^2\mf(st\n\phi_1+s(1-t)\n\phi_2+s\n\varphi_0
+\n\psi_0)
(t\p_i\phi_1+(1-t)\p_i\phi_2+\n\varphi_0)\mathrm{d}s\mathrm{d}t\\
&\quad\times(\p_j\phi_1-\p_j\phi_2).
\end{aligned}
\end{equation*}
This implies
\begin{equation}\label{2-39}
|I[\phi_1;\varphi_0,\psi_0]-I[\phi_2;\varphi_0,\psi_0]|
\leq C(1+\|\phi_1\|_{\ml}+\|\phi_2\|_{\ml})
(\|\phi_1-\phi_2\|_{\ml}).
\end{equation}
 Since $ \psi_0=\psi_{10}+\psi_{20} $, $ \varphi_0=-\zeta\psi_{20} $, so the continuity of $ I[\phi;\varphi_0,\psi_0] $ on $ \psi_{10} $ and $ \psi_{20} $ follows from the equality below
\begin{equation}\label{2-40}
\begin{aligned}
I[\phi;\varphi_{0},\psi_{0}]&
=\int_{0}^{1}t\p_{p_ip_j}^2\mf(\n\psi_{0}+(1-t)(\n\phi+\n\varphi_{0})
\p_i(\phi+\varphi_{0}) \p_j(\phi+\varphi_{0})\mathrm{d}t\\
&=\int_{0}^{1}t\p_{p_ip_j}^2\mf(\n\psi_{10}+\n\psi_{20}
+(1-t)(\n\phi-\n(\zeta\psi_{20})))\p_i(\phi-\zeta\psi_{20}) \p_j(\phi-\zeta\psi_{20})\mathrm{d}t.\\
\end{aligned}
\end{equation}
 \par Now, for each fixed $ \kappa_1 $, we will prove the continuous of the solution on  $ \kappa_2 $.  Let $ \kappa_2^m $ be a convergent sequence such that $ \kappa_2^m\rightarrow \kappa_2^\ast $ as $ m\rightarrow\infty $. Denote  $ \psi^m=\phi^m+\varphi_0^m+\psi_{10}+\psi_{20}^m=\phi^m+\psi_{10}+(1-\zeta)\psi_{20}^m $, where $\psi_{20}^m=-\kappa_2^m\int_{1}^{r}\rho_b(s)\frac{1}{s}\mathrm{d}s $. First we show that $ \phi^m\rightharpoonup \phi_{2}^\ast $ in $ \ml $ by using that $ \phi^m $ is a minimizing sequence of $ I[\phi;-\zeta\psi_{20}^\ast,\psi_{10}+\psi_{20}^\ast] $.
Since $ \kappa_2^m\rightarrow \kappa_{2}^\ast $ as $ m\rightarrow \infty $, so $ (1-\zeta)\psi_{20}^m\rightarrow (1-\zeta)\psi_{20}^\ast $, where $ \psi_{20}^\ast=-\kappa_2^\ast\int_{1}^{r}\rho_b(s)\frac{1}{s}\mathrm{d}s $.  All estimates in Lemma 2.3 and 2.4 as well as \eqref{2-9}  can be taken uniformly. In particular, $ \int_{\Omega} |\nabla\phi^m|^2\mathrm{d}x $ and $ \max_{\Omega}|\nabla\phi^m| $ are uniformly bounded.
\par For any given $ \delta>0 $ and sufficiently large $ m $, it follows from \eqref{2-39} and \eqref{2-40} that we have
\begin{equation*}
\left\{\begin{array}{l}
{\begin{aligned}
&|I[\phi^m;-\zeta\psi_{20}^m,\psi_{10}+\psi_{20}^m]
-I[\phi^m;-\zeta\psi_{20}^\ast,\psi_{10}+\psi_{20}^\ast]|<\delta,\\
&|I[\phi_{2}^\ast;-\zeta\psi_{20}^m,\psi_{10}+\psi_{20}^m]
-I[\phi_{2}^\ast;-\zeta\psi_{20}^\ast,\psi_{10}+\psi_{20}^\ast]|
<\delta.
\end{aligned}}
\end{array}\right.
\end{equation*}
Combing with the minimality of $ \phi^m $ for $ I[\phi;-\zeta\psi_{20}^m,\psi_{10}+\psi_{20}^m] $, we derive that
\begin{equation*}
I[\phi^m;-\zeta\psi_{20}^\ast,\psi_{10}+\psi_{20}^\ast]\leq
I[{\phi}^m;-\zeta\psi_{20}^m,\psi_{10}+\psi_{20}^m] +\delta
\leq I[\phi_{2}^\ast;-\zeta\psi_{20}^\ast,\psi_{10}+\psi_{20}^\ast]+2\delta.
\end{equation*}
Therefore, $ \phi^m $ is a minimizing sequence of $ I[\phi;-\zeta\psi_{20}^\ast,\psi_{10}+\psi_{20}^\ast] $. By the proof of existence of minimizer, we know that $ \phi^m\rightharpoonup \phi_{2}^\ast $ in $ \ml $. The uniformly convergence of $ \nabla \psi^m $ to $ \nabla \psi_{2}^\ast $ follows from uniform estimates in Lemma 2.3, 2.4  and a contradiction argument. Then we conclude that $ \max_{\Omega}|\n\psi^m|\rightarrow \max_{\Omega}|\n\psi_{2}^\ast| $. Hence $ \max_{\Omega}|\n\psi| $ is a continuous function of $ |\kappa_2| $.
\par  Finally, for each fixed $ \kappa_2 $, we will prove the continuous of the solution on  $ \kappa_1 $. Let $ \kappa_1^n $ be a convergent sequence such that $ \kappa_1^n\rightarrow \kappa_{1}^\star $ as $ n\rightarrow\infty $.
  Denote  $ \psi^n=\phi^n+\varphi_0+\psi_{10}^n+\psi_{20} $. Similar to the above proof, we have $ \max_{\Omega}|\n\psi^n|\rightarrow \max_{\Omega}|\n\psi_{1}^\star| $. Hence $ \max_{\Omega}|\n\psi| $ is a continuous function of $ |\kappa_1| $.
 \subsection{Removal of the truncation}\noindent
 \par In this subsection, we remove the truncation and complete the proof of Theorem 1.1 and Theorem 1.3. We only need prove Theorem 1.1, the proof of Theorem 1.3 is similar.
 \par Up to now, we have shown for fixed parameter $ \varepsilon $, \eqref{2-2} has a unique classical solution $ \psi $.  Let $\{\varepsilon_i\}_{i=1}^{\infty} $ be a strictly decreasing sequence such that $ \varepsilon_i\rightarrow 0 $ as $ i\rightarrow \infty $. For fixed $ \kappa_1 $ and $ i $, there exists a maximum interval $ [0, \kappa_2^i) $ such that for $ |\kappa_2| \in [0, \kappa_2^i) $,   $ |\n \psi|<\sqrt{1-2\varepsilon_i} $, Then   $ \psi $ is the solution  to the origin equation \eqref{1-18}. From the uniqueness of \eqref{2-2}, we can see $ \kappa_2^{i}\leq \kappa_2^{j} $ for $ i\leq j $. So $ \{\kappa_2^i \}_{i=1}^{\infty} $ is an increasing sequence with the upper bound $ 1 $, which implies the convergence of the sequence. Therefore, we  have
 \begin{equation}\label{2-41}
\hat \kappa_2=\lim_{i\rightarrow \infty}\kappa_2^{i}.
\end{equation}
If $ \kappa_2^i< \hat\kappa_2  $ for any $ i $, then for any $ |\kappa_2| \in [0,\hat \kappa_2 )$, there exists an index $ i $ such that $ |\kappa_2| \in [0,\kappa_2^i ) $, so the truncation can be removed such that $ \psi  $ is  the classical solution of \eqref{1-18} and $ \sup_{x\in\Omega}|\n \psi|<1 $.
\par The uniqueness of  subsonic spiral flows is easy to obtain. Suppose  there are two solutions to \eqref{1-18} such that $ \psi_i=\phi_i+\varphi_0+\psi_0 $, $ i=1,2 $ with $ \phi\in\ml $ and $ \max_{\Omega}|\n \psi_i |< 1 $. A small $ \varepsilon $ can be picked such that $ \max_{\Omega}|\n\psi |<\sqrt{1-2\varepsilon}  $. Both solutions will be solutions to the modified equation \eqref{2-2} with that $ \varepsilon $. It follows from the uniqueness of \eqref{2-2} that $ \psi_1=\psi_2 $. Therefore  the proof of  Theorem 1.1 is completed.
\par Moreover, we have $ \max_{\Omega}|\n\psi |\rightarrow 1$ as $ |\kappa_2|\rightarrow \hat \kappa_2 $. It is expected that subsonic spiral flows will tend to some  subsonic-sonic spiral flows as $ |\kappa_2|\rightarrow \hat \kappa_2 $. we will study this limiting behavior by compensated compactness framework.
\section{Subsonic-sonic spiral flows }\noindent

\par In this section, similar to the subsonic case, we only need prove Theorem 1.2.
 Firstly, let us recall the compensated compactness framework for steady irrotational flows in \cite{HWW11}.
\begin{theorem}
Let  $ \mathbf{u^\epsilon}(x_1,x_2)=(u_1^\epsilon,u_2^\epsilon)(x_1,x_2) $ be sequence of functions satisfying the following set of conditions $ (A) $:
 \item [$\rm{(A.1)}$]
 $ q^\epsilon(x_1,x_2)=|\mathbf{u^\epsilon}(x_1,x_2)|\leq 1 $ a.e. in $ \Omega $.
 \item [$\rm{(A.2)}$]
 $ \curl \mathbf{u^\epsilon} $ and $ \div (g(q^\epsilon)^2)\mathbf{u^\epsilon}) $ are confined in a compact set in $ H^{-1}_{loc}(\Omega) $.
 \\ Then there exists a subsequence (still labeled) $ \mathbf{u^\epsilon} $ that converges a.e. as $ \epsilon\rightarrow 0$ to $ \mathbf{\hat u} $ satisfying
  \begin{equation*}
 \hat q(x_1,x_2)=|\mathbf{\hat u}(x_1,x_2)|\leq 1 \quad a.e. \ (x_1,x_2)\in \Omega.
  \end{equation*}
 \end{theorem}
 \par Let $ (\rho^\epsilon,u_1^\epsilon,u_2^\epsilon )$ denote the solutions obtained in Theorem 1.1 to problem $ (\kappa_1,\kappa_2^\epsilon) $.
 Then we have
 \begin{equation}\label{3-1}
 \left\{\begin{array}{l}
\p_{x_1}((g(q^\epsilon)^2) u_1^\epsilon)+\p_{x_2}((g(q^\epsilon)^2) u_2^\epsilon)=0,\\
\p_{x_2}u_1^\epsilon-\p_{x_1}u_2^\epsilon=0.\\
\end{array} \right.
\end{equation}
  Thus conditions  $(A1) $ and $ (A2) $ in Theorem 3.1 are all satisfied. Theorem 3.1 implies that the solution sequence has a subsequence (still denoted by) $ (\rho^\epsilon,u_1^\epsilon,u_2^\epsilon )$ that converges a.e.  to a vector function $ (\hat\rho,\hat u_1,\hat u_2) $. Then the boundary condition
are satisfied for $ \hat\rho\hat{\mathbf{u}}$ in the sense of Chen-Frid \cite{CF99}.
  Since $\eqref{1-1}_2 $ and $\eqref{1-1}_3 $ hold for the sequence of subsonic solutions $ (\rho^\epsilon,u_1^\epsilon,u_2^\epsilon )$, so it is easy to see that $ (\hat \rho, \hat u_1, \hat u_2) $ also satisfies $\eqref{1-1}_2 $ and $\eqref{1-1}_3 $ in the sense of distribution. Thus the proof of Theorem 1.2 is completed.
 \par {\bf Acknowledgement.} Weng is partially supported by National Natural Science Foundation of China 11701431, 11971307, 12071359.

\end{document}